\documentclass [12pt]{article}
\usepackage{graphicx,amssymb,amsfonts,latexsym,amsmath,amsthm,times}
\usepackage{epsfig}
\usepackage{color}
\usepackage[english]{babel}
\usepackage[a4paper]{geometry}
\def\lanbox{\hbox{$\, \vrule height 0.25cm width 0.25cm depth 0.01cm \,$}}

\numberwithin{equation}{section}

\begin{document}

\vspace*{1.4cm}

\normalsize \centerline{\Large \bf THE PRESERVATION OF NONNEGATIVITY
OF SOLUTIONS}

\medskip

\centerline{\Large\bf OF A PARABOLIC SYSTEM WITH THE CUBED LAPLACIAN}

\vspace*{1cm}

\centerline{\bf Messoud Efendiev$^{1,2}$, Vitali Vougalter$^{3 \ *}$}

\bigskip

\centerline{$^1$ Helmholtz Zentrum M\"unchen, Institut f\"ur Computational
Biology, Ingolst\"adter Landstrasse 1}

\centerline{Neuherberg, 85764, Germany}

\centerline{e-mail: messoud.efendiyev@helmholtz-muenchen.de}

\centerline{$^2$ Azerbaijan University of Architecture and Construction, Baku, Azerbaijan}

\centerline{e-mail: messoud.efendiyev@gmail.com}

\bigskip

\centerline{$^{3 \ *}$  Department of Mathematics, University
of Toronto}

\centerline{Toronto, Ontario, M5S 2E4, Canada}

\centerline{ e-mail: vitali@math.toronto.edu}

\medskip


\vspace*{0.25cm}

\noindent {\bf Abstract:}
The article is devoted to the easily verifiable necessary condition of the
preservation of the nonnegativity of the solutions of a system of parabolic
equations containing the cubed Laplacian. Such necessary condition
is extremely important for the applied analysis community since it imposes
the necessary form of the system of equations that must be studied
mathematically.

\vspace*{0.25cm}



\vspace*{0.5cm}

\bigskip

\bigskip


\setcounter{section}{1}

\centerline{\bf 1. Introduction}

\medskip

The solutions of many systems of convection-diffusion-reaction problems
appearing in biology, physics or engineering describe such quantities as
population densities, pressure or concentrations of nutrients and chemicals.
Therefore, a natural property to require for the solutions is their nonnegativity.
The models which fail to guarantee the nonnegativity are not valid or break down
for small values of the solution. In many situations, showing that a particular
model does not preserve the nonnegativity leads to the better understanding
of the model and its limitations. One of the first steps in analyzing
ecological or biological or bio-medical models mathematically is to test
whether the solutions which originate from the nonnegative initial data remain
nonnegative (as long as they exist). In other words, the model under
consideration ensures that the nonnegative cone is positively invariant.
Let us recall that if the solutions (of a given evolution PDE) which correspond
to the nonnegative initial data remain nonnegative as long as they exist,
we say that such system satisfies the nonnegativity property.

For scalar equations involving the standard Laplace operator the nonnegativity
property follows directly from the maximum principle (see ~\cite{E13} and
the references therein). However, for the equations containing the cubed Laplacian
the maximum principle does not hold. 

In the present work we aim to establish a simple and easily verifiable criterion, that
is, the necessary condition for the nonnegativity of solutions
of systems of nonlinear convection--diffusion--reaction equations involving the
cubed Laplacian appearing in the modelling of life sciences. We believe that
it could provide the modeler with a tool, which is easy to verify, to
approach the problem of the nonnegative invariance of such model. 

Currently we study the preservation of the nonnegativity
of solutions of the following system of reaction-diffusion equations
\begin{equation}
\label{h}
\frac{\partial u}{\partial t} =
A {\Delta}^{3} u + \sum_{i=1}^{d}\Gamma^{i}
\frac{\partial u}{\partial x_{i}}-F(u), \quad x\in {\mathbb R}^{d}, \quad
d\in {\mathbb N},
\end{equation}
where $A, \ \Gamma^{i}, \ 1\leq i\leq d$ are $N\times N$ matrices with constant
coefficients. Here
$$
u(x,t)=(u_{1}(x,t), u_{2}(x,t), ..., u_{N}(x,t))^{T}.
$$
Note that in the present work we consider the space of an arbitrary
dimension. Solvability
conditions for a linearized Cahn-Hilliard equation of sixth order
were obtained in ~\cite{VV12}.
The solvability of the single equation containing the standard Laplacian with
transport important for the fluid mechanics was discussed in ~\cite{VV121}.
The verification of biomedical processes with anomalous diffusion, transport and interaction of species 
 in the situation of the one spatial dimension with a single Laplacian raised to the power  $\displaystyle{0 < s < \frac{1}{4}}$
 in the diffusion term was accomplished in ~\cite{EV21} (see also ~\cite{EV22}).
The similar ideas in the space of dimension
$d\geq 2$ in the case when the diffusion term of the system
 involves the sum of the standard Laplace operator, which acts on the first $m$ variables and the
 fractional Laplacian with respect to the remaining  variables were used in ~\cite{EV211}.
The necessary conditions for preserving the nonnegative cone in the case of the double scale anomalous diffusion
were derived in ~\cite{EV221}. 
The necessary and sufficient conditions for an infinite system of parabolic equations preserving the positive cone were discussed
in ~\cite{ELS10}. The positivity of solutions of systems of stochastic PDEs was covered in ~\cite{CES13}. Article ~\cite{V22} is devoted to
the  preservation of nonnegativity of solutions of a parabolic system with the bi-Laplacian.
 In ~\cite{EP07} the
authors study the large time behavior of solutions of a class of fourth order parabolic equations defined on unbounded domains
using the Kolmogorov $\varepsilon$-entropy as a measure. Spatial patterns arising in higher order models in physics and mechanics
were treated in ~\cite{PT01}.
Let us  assume here that system (\ref{h}) contains the square matrices with 
the entries constant in space and time
$$
(A)_{k,j}:=a_{k,j}, \quad ({\Gamma}^{i})_{k,j}:=\gamma^{i}_{k,j}, \quad 1\leq k,j
\leq N, \quad 1\leq i\leq d, \quad d\in {\mathbb N}
$$
and that the matrix $A+A^{*}>0$ for the sake of the global well posedness
of  (\ref{h}). Here $A^{*}$ denotes the adjoint of matrix $A$.
Thus, the system of equations (\ref{h}) can be rewritten in the form
\begin{equation}
\label{h1}  
\frac{\partial u_{k}}{\partial t}=\sum_{j=1}^{N}a_{k,j}{\Delta}^{3}u_{j}+
\sum_{i=1}^{d}\sum_{j=1}^{N}\gamma^{i}_{k,j}\frac{\partial u_{j}}
{\partial x_{i}}-F_{k}(u), \quad 1\leq k\leq N.
\end{equation}
In the present article the interaction vector function term
$$
F(u):=(F_{1}(u), F_{2}(u), ..., F_{N}(u))^{T},
$$
which can be linear or nonlinear. We assume
its smoothness in the theorem further down for the sake of the well posedness of
our system of equations (\ref{h}), although, we are not focused on the well posedness
issue in the current work.
Let us introduce the inner product as
\begin{equation}
\label{ip}
(f(x), g(x))_{L^{2} ({\mathbb R}^{d})}:=\int_{{\mathbb R}^{d}}f(x){\bar g(x)}dx.
\end{equation}
As for the
vector functions, their inner product is defined using their components as
\begin{equation}
\label{ipv}  
(u, v)_{L^{2} ({\mathbb R}^{d}, {\mathbb R}^{N})}:=\sum_{k=1}^{N}(u_{k}, v_{k})_
{L^{2}({\mathbb R}^{d})}.
\end{equation}
Obviously, (\ref{ipv}) induces the norm
$$
\|u\|_{L^{2} ({\mathbb R}^{d}, {\mathbb R}^{N})}^{2}=\sum_{k=1}^{N}\|u_{k}\|_
{L^{2}({\mathbb R}^{d})}^{2}.
$$
By the nonnegativity of a vector function below we mean the nonnegativity
of the each of its components. Our main result is as
follows.

\bigskip

\noindent
{\bf Theorem 1.} {\it Let $F: {\mathbb R}^{N}\to {\mathbb R}^{N}$,
so that $F\in {\mathbb C}^{1}$, the initial condition for system
(\ref{h}) is
$u(x,0)=u_{0}(x)\geq 0$ and $u_{0}(x)\in L^{2} ({\mathbb R}^{d}, {\mathbb R}^{N})$.
Let us also assume that the off diagonal elements of the matrix $A$, are
nonnegative, so that
\begin{equation}
\label{akl}
a_{k, j}\geq 0, \quad  1\leq k,j\leq N, \quad k\neq j.
\end{equation}   
Then the necessary condition for the system of equations (\ref{h}) to admit a solution
$u(x,t)\geq 0$ for all $t\in [0, \infty)$ is that the matrices $A$ and
$\Gamma^{i}, \ 1\leq i\leq d$ are diagonal and for all $1\leq k\leq N$
\begin{equation}
\label{fc}  
F_{k}(s_{1}, ..., s_{k-1}, 0, s_{k+1}, ..., s_{N})\leq 0
\end{equation}  
is valid, where $s_{l}\geq 0$  
and} $1\leq l\leq N, \ l\neq k$.

\bigskip

\noindent
{\bf Remark 1.} {\it In the situation of the linear interaction term, namely
when $F(u)=Lu$, where $L$ is a matrix with elements
$b_{i,j}, \ 1\leq i,j\leq N$ constant in space and time, our necessary condition
yields the condition that the matrix $L$ must be essentially nonpositive,
that is $b_{i,j}\leq 0$ for  $i\neq j, \ 1\leq i,j\leq N$.}

\bigskip

\noindent
{\bf Remark 2.} {\it Our proof implies that, the necessary condition for
preserving the nonnegative cone is carried over from the ODE (the spatially
homogeneous case, as described by the ordinary differential equation
$u'(t)=-F(u)$) to the situation of the diffusion containing the cubed Laplacian and the
convective transport term.}

\bigskip

\noindent
{\bf Remark 3.} {\it In the consecutive articles we plan to discuss the
following situations:

\noindent
a) the necessary and sufficient conditions of the current work,

\noindent
b) the nonautonomous version of the present article,

\noindent 
c) the density-dependent diffusion matrix,

\noindent 
d) the effect of the delay term in the situations a), b) and c).}
      
\bigskip

Let us turn our attention to the proof of our main proposition.

\bigskip


\setcounter{section}{2}

\centerline{\bf 2. The preservation of the nonnegativity of the solution
of the system of parabolic equations}

\bigskip

\noindent
{\it Proof of Theorem 1.} We note that the maximum principle broadly
used to study the solutions of single parabolic equations with the
standard Laplace operator is not applicable to these problems involving the cubed
Laplacian. Let us consider a time independent, square integrable, nonnegative
vector function $v(x)$ and estimate
$$
\Bigg(\frac{\partial u}{\partial t}\Bigg|_{t=0}, v \Bigg)_
{L^{2} ({\mathbb R}^{d}, {\mathbb R}^{N})}
=\Bigg(\hbox{lim}_{t\to 0^{+}}\frac{u(x,t)-u_{0}(x)}{t}, v(x)\Bigg)_
{L^{2} ({\mathbb R}^{d}, {\mathbb R}^{N})}.
$$
By means of the continuity of the inner product, the right side of the
identity above is equal to
\begin{equation}
\label{uv}
\hbox{lim}_{t\to 0^{+}}\frac{(u(x,t),v(x))_{L^{2} ({\mathbb R}^{d}, {\mathbb R}^{N})}}{t}-
\hbox{lim}_{t\to 0^{+}}\frac{(u_{0}(x),v(x))_{L^{2} ({\mathbb R}^{d}, {\mathbb R}^{N})}}{t}.
\end{equation}
We choose the initial condition for our system of equations $u_{0}(x)\geq 0$ and the
constant in time vector function $v(x)\geq 0$ to be orthogonal to each other in 
$L^{2} ({\mathbb R}^{d}, {\mathbb R}^{N})$. It can be accomplished, for example for
\begin{equation}
\label{u0vj}  
u_{0}(x)=({\tilde u_{1}}(x), ..., {\tilde u_{k-1}}(x), 0,
{\tilde u_{k+1}}(x), ..., {\tilde u_{N}}(x)), \quad v_{j}(x)={\tilde v}(x)
\delta_{j, k}, 
\end{equation}
with $1\leq j\leq N$. Here $\delta_{j, k}$ stands for the Kronecker symbol and
$1\leq k\leq N$ is fixed. Therefore, the second term in (\ref{uv})
is trivial and (\ref{uv}) equals to
$$
\hbox{lim}_{t\to 0^{+}}\frac{\sum_{k=1}^{N}\int_{{\mathbb R}^{d}}u_{k}(x,t)v_{k}(x)dx}
{t}\geq 0,
$$
because all the components $u_{k}(x,t)$ and $v_{k}(x)$ contained in the formula
above are nonnegative. Thus, we arrive at
$$
\sum_{j=1}^{N}\int_{{\mathbb R}^{d}}\frac{\partial u_{j}}{\partial t}\Bigg|_{t=0}
v_{j}(x)dx\geq 0.
$$
By virtue of (\ref{u0vj}), only the $k$ th component of the vector function
$v(x)$ does not vanish identically. This yields
$$
\int_{{\mathbb R}^{d}}\frac{\partial u_{k}}{\partial t}\Bigg|_{t=0}{\tilde v}(x)dx
\geq 0.
$$
Let us use (\ref{h1}) to derive
$$
\int_{{\mathbb R}^{d}}\Bigg[\sum_{j=1, \ j\neq k}^{N}a_{k, j}{\Delta}^{3}{\tilde u_{j}}(x)+
\sum_{i=1}^{d}\sum_{j=1, \ j\neq k}^{N}\gamma^{i}_{k, j}\frac{\partial {\tilde u_{j}}}
{\partial x_{i}}-
$$
$$
F_{k}({\tilde u_{1}}(x) , ..., {\tilde u_{k-1}}(x), 0,
{\tilde u_{k+1}}(x), ...,  {\tilde u_{N}}(x))\Bigg]{\tilde v}(x)dx
\geq 0.
$$
Because the nonnegative, square integrable function ${\tilde v}(x)$ can be
chosen arbitrarily, we obtain
$$
\sum_{j=1, \ j\neq k}^{N}a_{k, j}{\Delta}^{3}{\tilde u_{j}}(x)+\sum_{i=1}^{d}
\sum_{j=1, \ j\neq k}^{N}{\gamma}^{i}_{k, j}\frac{\partial {\tilde u_{j}}}
{\partial x_{i}}-
$$
\begin{equation}
\label{sc0}  
F_{k}({\tilde u_{1}}(x) , ..., {\tilde u_{k-1}}(x), 0,
{\tilde u_{k+1}}(x), ...,  {\tilde u_{N}}(x))    \geq 0 \quad a.e.
\end{equation}
For the purpose of the scaling, we replace all the ${\tilde u_{j}}(x)$ by
$\displaystyle{{\tilde u_{j}}\Bigg(\frac{x}{\varepsilon}\Bigg)}$ in the
bound above, where $\varepsilon>0$ is a small parameter. This implies that
$$
\sum_{j=1, \ j\neq k}^{N}\frac{a_{k, j}}{{\varepsilon}^{6}}{\Delta}^{3}{\tilde u_{j}}
(y)+\sum_{i=1}^{d}\sum_{j=1, \ j\neq k}^{N}\frac{\gamma^{i}_{k, j}}{\varepsilon}\frac
{\partial {\tilde u_{j}}(y)}{\partial y_{i}}-
$$
\begin{equation}
\label{sc}  
F_{k}({\tilde u_{1}}(y) , ..., {\tilde u_{k-1}}(y), 0,
{\tilde u_{k+1}}(y), ...,  {\tilde u_{N}}(y))    \geq 0 \quad a.e.
\end{equation}
Let us first suppose that some of the $a_{k, j}$ contained in the sum in the left side
of (\ref{sc}) are strictly positive.
Clearly, the first term in the left side of (\ref{sc}) is the leading one
as $\varepsilon\to 0$. We choose here all the
${\tilde u_{j}}(y), \ 1\leq j\leq N, \ j\neq k$ to be identical, equal to
$\displaystyle{2-e^{\sum_{k=1}^{d}y_{k}}}$ in a neighborhood of the origin, smooth and decaying to zero at the
infinities. A trivial calculation gives us
$$
{\Delta}^{3}{\tilde u_{j}}(y)\Big|_{y=0}=-d^{3}<0, \quad 1\leq j\leq N, \quad
j\neq k.
$$
This means that ${\Delta}^{3}{\tilde u_{j}}(y)<0$ in a neighborhood of the origin
due to the simple continuity argument.
By making the parameter $\varepsilon$ sufficiently small, we are able to violate 
inequality in (\ref{sc}). Since the negativity of the off diagonal elements
of the matrix $A$ is ruled out by means of assumption (\ref{akl}), we have
$$
a_{k, j}=0, \quad 1\leq k, j\leq N, \quad k\neq j.
$$
Thus, from (\ref{sc}) we derive
$$ 
\sum_{i=1}^{d}\sum_{j=1, \ j\neq k}^{N}\frac{\gamma^{i}_{k, j}}{\varepsilon}
\frac{\partial}{\partial y_{i}}{\tilde u_{j}}(y)-
$$
\begin{equation}
\label{sc1} 
F_{k}({\tilde u_{1}}(y) , ..., {\tilde u_{k-1}}(y), 0,
{\tilde u_{k+1}}(y), ...,  {\tilde u_{N}}(y))    \geq 0 \quad a.e.
\end{equation}
We choose here
$$
{\tilde u_{j}}(y)=Q_{j, i}e^{-y_{i}sign \gamma^{i}_{k, j}}
$$
 in a neighborhood of the origin, smooth and
decaying to zero at the infinities with $Q_{j, i}(y)$ positive and independent of $y_{i}$.
Then the left side of (\ref{sc1}) can be made as negative as possible when
$\varepsilon\to 0$.
This will violate inequality (\ref{sc1}). Evidently,
the last term in the left side of (\ref{sc1}) will remain bounded. 
This means that for the matrices $\Gamma^{i}$ involved in the system of equations (\ref{h}), the off
diagonal elements should vanish, such that
$$
\gamma^{i}_{k, j}=0, \quad 1\leq k,j\leq N, \quad k\neq j, \quad 1\leq i\leq d.
$$
Therefore, by virtue of (\ref{sc1}) we obtain that
$$
F_{k}({\tilde u_{1}}(x), ..., {\tilde u_{k-1}}(x), 0, {\tilde u_{k+1}}(x), ...,
{\tilde u_{N}}(x))\leq 0 \quad a.e.
$$
with ${\tilde u_{j}}(x)\geq 0$ and ${\tilde u_{j}}(x)\in L^{2}({\mathbb R}^{d})$,
where $1\leq j\leq N, \ j\neq k$.   \hfill\lanbox

\bigskip

\noindent
{\bf Acknowledgement.} {\it  V.V. is grateful to Israel Michael Sigal for the partial support by the NSERC grant
 NA 7901.}

\bigskip


\end{document}